\documentclass[11pt]{article}
\usepackage{amscd,amsfonts,amsmath,amsthm,euscript}

\def\p{\partial}
\def\RE{\mathbb{R}}
\def\CX{\mathbb{C}}
\def\ZE{\mathbb{Z}}
\def\HH{\mathcal{H}}
\def\bfn{\mathbf{n}}
\def\bfv{\mathbf{v}}
\def\eps{\varepsilon}
\def\we{\wedge}
\def\AM{\EuScript{A}}

\DeclareMathOperator{\vol}{vol}
\DeclareMathOperator{\Vol}{Vol}
\DeclareMathOperator{\Ind}{index}
\DeclareMathOperator{\Spin}{Spin}

\renewcommand{\Re}{\operatorname{Re}}
\renewcommand{\Im}{\operatorname{Im}}

\newtheorem{thm}{Theorem}
\newtheorem{cor}[thm]{Corollary}
\newtheorem{prop}[thm]{Proposition}

\theoremstyle{definition}
\newtheorem{defi}[thm]{Definition}
\newtheorem{example}[thm]{Example}
\theoremstyle{remark}
\newtheorem*{remark}{Remark}
\newtheorem*{remarks}{Remarks}

\begin{document}

\title{Deformations of calibrated submanifolds with boundary}
\author{{\sc Alexei Kovalev}\\[5pt]
DPMMS, University of Cambridge}
\date{}
\maketitle

\begin{abstract}
We review some results concerning the deformations of calibrated
minimal submanifolds which occur in Riemannian manifolds with special
holonomy. The calibrated submanifolds are assumed compact with a non-empty
boundary which is constrained to move in a particular fixed submanifold. 
The results extend McLean's deformation theory previously developed for
closed compact submanifolds.
\end{abstract}

\section{Preliminaries}

Calibrated submanifolds are a particular type of minimal submanifolds
introduced by Harvey and Lawson~\cite{HL} as a generalization of complex
submanifolds of K\"ahler manifolds; for a detailed reference, see
{\it op.cit.} and \cite{harvey}. Harvey and Lawson \cite{HL} also found four
new types on calibrations defined on Euclidean spaces and, more generally, on
Ricci-flat Riemannian manifolds with reduced holonomy. Further examples of
calibrations were subsequently discovered, see \cite[\S 4.3]{joyce} and
references therein.

McLean~\cite{mclean} studied deformations for the four types of
calibrated submanifolds defined in~\cite{HL} and showed that the deformation
problem may be interpreted as a non-linear PDE with Fredholm properties.
In some cases, there is a smooth finite-dimensional `moduli space'.
The submanifolds in~\cite{mclean} were assumed compact and without boundary.
McLean's deformation theory in~\cite{mclean} was later extended by several
authors to more general classes of submanifolds. In this paper, we survey 
the generalizations of McLean's results to compact submanifolds having
non-empty boundary. Proofs will at most be briefly sketched or referred to
the original papers.

We begin in this section with some key concepts of calibrated geometry
and foundational results concerning deformations of compact submanifolds
(possibly with boundary).
The remainder of the paper is organized in four sections dealing with 
the four types of calibrated submanifolds in \cite{HL,mclean}, namely
special Lagrangian, coassociative, associative and Cayley submanifolds.

Submanifolds are taken to be embedded and connected, for notational
convenience (it can be checked that the results extend to immersed
submanifolds). Smooth functions and, more generally, sections of vector
bundles on (sub)manifolds with boundary are understood as `smooth up to
the boundary', so at each point of the boundary these have one-sided partial
derivatives of any order in the inward-pointing normal direction.

\subsection{Calibrations}

\begin{defi}
Let $(M,g)$ be a Riemannian manifold. For any tangent $k$-plane $V$,
i.e. a $k$-dimensional subspace $V$ of a tangent space $T_xM$, a choice of
orientation on~$V$ together with the restriction of $g$ determines a
natural {\em volume form} on $V$, $\vol_V\in \Lambda^kV^*$.

A differential $k$-form $\phi$ on $M$ is called a {\em calibration} if
(i) $d\phi = 0$ and 
(ii) for each $x\in M$ and every oriented $k$-dimensional subspace
$V\subset T_xM$, $\phi|_V=a \vol_V$ for some $a\le 1$.

An oriented  $k$-dimensional submanifold $N$ of $M$ is said to be
{\em calibrated by $\phi$}  if the pull-back of
$\phi$ to $N$ coincides with the Riemannian volume form for the metric on $N$
induced by $g$, i.e.\ $\phi|_{T_x N}=\vol_{T_x N}$ for each $x\in N$.
\end{defi}

The next result shows that calibrated submanifolds are minimal, in fact
volume-minimizing (if compact).

\begin{prop}\label{calibr-min}
Let $M$ be a Riemannian manifold and let $\phi\in\Omega^k(M)$ be a calibration
on~$M$.

\begin{itemize}
\item[\emph{(a)}]
If a closed $k$-dimensional submanifold $X\subset M$ is calibrated by
$\phi$, then $X$ is volume-minimizing in its homology class. Moreover, if $Y$
is a volume-minimizing closed $k$-dimensional submanifold of~$M$ in the
homology class of $X$, then $Y$ is calibrated by~$\phi$.
\item[\emph{(b)}]
Let $W\subset M$ be a submanifold such that $\phi|_W=0$. 
If $X\subset M$ is a calibrated compact $k$-dimensional submanifold with
non-empty boundary $\p X\subset W$, then $X$ is volume-minimizing in the
relative homology class $[X]\in H_k(M,W;\ZE)$. Moreover, if $Y$ is a
compact submanifold of~$M$ with boundary $\p Y\subset W$ and $Y$ is
volume-minimizing in the relative homology class $[X]$, then $Y$ is calibrated
by~$\phi$.
\end{itemize}
\end{prop}

The clause (a) is proved in \cite[Thm.~II.4.2]{HL} by application of
Stokes' theorem. The extension (b) to submanifolds with boundary follows by
a similar argument as the hypothesis $\phi|_W=0$ ensures the vanishing of
the additional terms arising from the boundary
(cf.~\cite[pp.~1233--1234]{gayet}).
Suppose that a submanifold $Y$ with boundary in $W$ represents the 
relative homology class $[X]$. Considering $X$ and $Y$ as chains, we can find
a $(k+1)$-dimensional chain $N$ with boundary in $W$ and a $k$-dimensional
chain $P$ contained in $W$ so that $Y-X=\p N + P$. We then obtain
$$
\Vol(Y) \ge \int_Y \phi = \int_X \phi + \int_{\p N}\phi + \int_P \phi = \Vol(X),
$$
noting in the inequality that $\phi$ is a calibration, then applying Stokes'
theorem and taking account of the vanishing of $\phi$ on~$W$.

The above argument may be viewed as a generalization of the volume-minimizing
property of compact complex submanifolds of K\"ahler manifolds by application
of Wirtinger's inequality. (In a K\"ahler manifold with K\"ahler form
$\omega$, every $k$-dimensional complex submanifold is calibrated by
$\omega^k/k!$.)

\subsection{Normal deformations of submanifolds}

McLean's deformation theory~\cite{mclean} was originally developed for closed
compact submanifolds. The nearby deformations of a given closed submanifold of
a Riemannian manifold may be assumed to be {\em normal deformations}, defined
using the Riemannian exponential map on normal vector fields, i.e.\ $C^1$
sections of the normal bundle of this submanifold.

When a submanifold has a boundary, we wish to consider the deformation problem
as an elliptic boundary value problem. We shall in fact require that the
boundary moves in a certain fixed submanifold which,
following~\cite{butscher}, we call a {\em scaffold}. In general, we cannot
use, as in the case of closed submanifolds, exponential deformations defined
using the given metric $g$ on the ambient manifold, since the scaffold may
not be preserved under such deformations. We shall define
on~$M$ a modified metric $\hat{g}$ whose associated exponential map does
preserve the scaffold because the scaffold will be totally geodesic with
respect to the new metric. (The actual construction of the metric $\hat{g}$
will depend on the considered calibration.)

\begin{prop}\label{tubular}
\mbox{}
\begin{itemize}
\item[\emph{(a)}] 
Let $P$ be a closed submanifold of a Riemannian manifold~$M$.

There exist
an open subset $V_P$ of the normal bundle $N_{P/M}$ of $P$ in~$M$,
containing the zero section, and a tubular neighbourhood $T_P$ of $P$
in~$M$, such that the exponential map $\exp_{M}|_{V_P}:V_P\to T_P$ is a
diffeomorphism onto~$T_P$.
\item[\emph{(b)}]
Let $M$ be a smooth manifold of dimension~$n$ and $P\subset M$ a compact
submanifold with non-empty boundary~$\p P$. Let $W$ be a submanifold of~$M$
with $\p P\subset W$ and let $\hat{g}$ be a Riemannian metric on~$M$ such that
$P$ and $W$ meet orthogonally and $W$ is totally geodesic with respect
to~$\hat{g}$.

There exists an open subset $V_P$ of the normal bundle $\hat{N}_{P/M}$ of
$P$ in~$M$, containing the zero section, and an $n$-dimensional submanifold
$T_P$ of $M$ with boundary such that $P\subset T_P$ and
$\widehat\exp_{M}|_{V_P}:V_P\to T_P$ is a diffeomorphism onto~$T_P$.
Furthermore, if a section $\bf{v}$ of $\hat{N}_{P/M}$ takes values in $V_P$,
then $\widehat\exp_{M}(\bfv(x))\in W$ for all $x\in\p P$.
\end{itemize}
\end{prop}

The clause (a) is a consequence of the tubular neighbourhood
theorem~\cite[Chapter IV, Thm.~9]{lang}. For the extension (b) to
submanifolds with boundary and the existence of the adapted metric $\hat{g}$
modifying a given metric on~$M$ on an open neighbourhood of $\p P$, cf.\
\cite[Prop. 6]{butscher} or \cite[Prop.~4.4]{KL}. 
We stress that the new auxiliary metric $\hat{g}$ in (b) is used
{\em solely} for the purpose of considering the exponential map and
applications of the tubular neighbourhood theorem --- but {\em not} for the
minimal or volume-minimizing properties of calibrated submanifolds (for
which we continue to use the original metric on~$M$).

By {\em nearby deformations} of a compact submanifold $P$ (with or without
boundary) we shall mean submanifolds of the form $P_{\bfv}=\exp_{\bf{v}}(P)$,
where $\bfv$ is a $C^1$-section of the normal bundle $N_{P/M}$. The
section $\bfv$ is assumed sufficiently small in the $C^1$ norm, so that
$P_{\bfv}$ is contained in a tubular neighbourhood $T_P$ defined by
Proposition~\ref{tubular}. We shall call sections of $N_{P/M}$ the
{\em normal vector fields} on~$P$.

We interpret submanifolds
as appropriate equivalence classes of the embedding (more generally,
immersion) maps and the term `moduli space of submanifolds' is used below in
this sense.

\section{Special Lagrangian submanifolds in Calabi--Yau manifolds}

Let $M$ be a K\"ahler manifold of complex dimension~$m$, with K\"ahler form
$\omega$ and suppose further that the metric on~$M$ has holonomy contained in
$SU(m)$. Then the canonical bundle of $M$ may be trivialized by a holomorphic
$(m,0)$-form $\Omega$ satisfying
\begin{equation}\label{omegas}
(-1)^{m(m-1)/2}\,(i/2)^m\,\Omega\we\bar\Omega=\omega^m/m! \; .
\end{equation}
It also follows that the K\"ahler metric is Ricci-flat. Conversely, if 
a K\"ahler form $\omega$ and a holomorphic $(m,0)$-form $\Omega$ satisfy
\eqref{omegas} on $M$, then this K\"ahler metric has holonomy in $SU(m)$.
We shall call $(M,\omega,\Omega)$ as above a {\em Calabi--Yau manifold}.

The real $m$-form $\Re\Omega$ is a calibration on~$M$
(cf. \cite[Thm. III.1.10]{HL}) and the submanifolds 
calibrated by $\Re\Omega$ are called {\em special Lagrangian submanifolds}.
It will be convenient to use an equivalent definition.
\begin{prop}[{\cite[Cor. III.1.11]{HL}}]\label{SLdefi}
Let $(M,\omega,\Omega)$ be a Calabi--Yau manifold of complex dimension~$m$.
A real $m$-dimensional submanifold $L\subset M$ with some choice of
orientation is special Lagrangian if and only if
\begin{equation}\label{SLeq}
\omega|_L=0 \qquad\text{and}\qquad \Im\Omega|_L=0.
\end{equation}
\end{prop}

The following result about the deformations of compact special Lagrangian
submanifolds without boundary is due to McLean.
\begin{thm}[{\cite[Thm.~3.6]{mclean}}]\label{mclean-SL}
Let $M$ be a Calabi--Yau manifold and $L$ a closed special Lagrangian
submanifold in~$M$. Then the moduli space of nearby special Lagrangian
deformations of~$L$ is a smooth manifold of dimension the first Betti
number $b^1(L)$.
\end{thm}

The argument of Theorem~\ref{mclean-SL} uses the equivalent
definition~\eqref{SLeq} of special Lagrangian submanifolds in terms of the
vanishing of differential forms.

Applying Proposition~\ref{tubular}(a) to $L$ we may write any nearby
deformation of $L$ as $L_v=\exp_v L$ for a section $v\in\Gamma (N_{L/M})$
of the normal bundle. On the other hand, there is
an isometry of vector bundles
\begin{equation}\label{1forms}
\jmath_L:\mathbf{v}\in N_{L/M}\to
(\mathbf{v}\,\lrcorner\,\omega)|_L\in\Lambda^1T^*L
\end{equation}
defined using the K\"ahler form on~$M$.
Thus the nearby deformations of $L$ are equivalently given by
`small' 1-forms on~$L$.

The map
\begin{equation}\label{SLmap}
F:\alpha\in\Omega^1(L)\to  (\exp_{\bfv}^*(\omega),\exp_{\bfv}^*(\Im\Omega))
\in\Omega^2(L)\oplus\Omega^m(L),
\quad  \bfv=\jmath_L^{-1}(\alpha),
\end{equation}
is defined for `small' $\alpha$, and $F(\alpha)=0$ precisely if
$\exp_{\bf{v}}(P)$ is a special Lagrangian deformation.

\begin{prop}\label{SL}
Let $L$ be a special Lagrangian submanifold of a Calabi--Yau manifold and
$F$ the `deformation map' defined in~\eqref{SLmap}.
\begin{itemize}
\item[\emph{(a)}]
The map $F$ is smooth, with derivative at $\alpha=0$ given by
$$
dF|_0(\alpha)=(d\alpha,d*\alpha),\qquad \alpha\in\Omega^1(L).
$$
\item[\emph{(b)}]
If $L$ is a closed submanifold, then there is a neighbourhood $T_L$ of the
zero 1-form such that the image $F(T_L)$ consists of pairs of exact forms.
\end{itemize}
\end{prop}
Note that in Proposition~\ref{SL}(a) a special Lagrangian $L$ need not be
compact.

The nearby special Lagrangian deformations of~$L$ correspond to the
1-forms $\alpha$ satisfying a non-linear differential equation
$F(\alpha)=0$ of first order, with $F(0)=0$.
One can show using Hodge theory for a closed manifold $L$ and the implicit
function theorem in Banach spaces that the $C^1$-small solutions 
of the special Lagrangian deformation problem may be parameterised by the
closed and co-closed 1-forms $\alpha$ on~$L$. By Hodge theory again, as $L$ is
a closed manifold, the moduli space of special Lagrangian deformations is
then locally parameterised by the vector space of harmonic 1-forms on $L$
or, equivalently, by the de Rham cohomology group $H^1(L,\RE)$ of dimension
$b^1(L)$. Theorem~\ref{mclean-SL} follows.

The constraint~\eqref{omegas} determines a holomorphic form $\Omega$ up to a
factor $e^{i\theta}$ for some real constant~$\theta$ and $\Re(e^{i\theta}\Omega)$
is also a calibration on~$M$. Manifolds calibrated by
$\Re(e^{i\theta}\Omega)$ are called special Lagrangian with phase $\theta$.
Every submanifold calibrated by $\Re(e^{i\theta}\Omega)$ is Lagrangian
(with respect to the symplectic form $\omega$) and minimal (with respect to
the metric on~$M$). 
Conversely, it is known that every connected minimal Lagrangian
submanifold in a Calabi-Yau manifold is calibrated by $\Re(e^{i\theta}\Omega)$
for some real constant~$\theta$ \cite[cf.\ Prop.~III.2.17]{HL}.

If $L'$ is a minimal Lagrangian deformation of $L$ then by the above
$L'$ must be calibrated and volume-mininizing. Therefore, $L$ is special
Lagrangian by application of Proposition~\ref{calibr-min}(a). It follows
that moduli space in Theorem~\ref{mclean-SL} can be  equivalently regarded
as the space of nearby minimal Lagrangian deformations of~$L$.
In particular, there is no loss of generality in restricting attention to
submanifolds calibrated by $\Re\Omega$, i.e. with $\theta=0$.

\begin{remark}
Salur~\cite{salur} extended the result of Theorem~\ref{mclean-SL} to the
situation when the almost complex structure on $M$ is not necessarily
integrable. More explicitly, $M$ is a Hermitian symplectic $2m$-manifold
with symplectic form $\omega$, the metric on $M$ is Hermitian with respect
to an $\omega$-compatible almost-complex structure and there is a
(non-vanishing) complex $(m,0)$-form $\Omega$ on~$M$ satisfying
\eqref{omegas}. The main theorem in~\cite{salur} then asserts that the
moduli space of nearby special Lagrangian deformations of $L$ with
arbitrary phase is smooth with dimension $b^1(L)$.
In the case when $d\Im\Omega=0$, it can be checked that any special Lagrangian
deformations of $L$ necessarily have the same phase as~$L$ (and can be
obtained essentially by McLean's argument in~\cite{mclean}).
\end{remark}

Suppose now that a compact special Lagrangian submanifold $L$ has non-empty
boundary $\p L$.
We shall need the following definition from~\cite[p.~1954]{butscher}.
\begin{defi}\label{SL-scaf}
Let $M$ be a Calabi--Yau manifold and
let $L\subset M$ be a submanifold with boundary $\p L$. Denote by
$\bfn\in\Gamma(T_{\p L} L)$ the inward unit normal vector field.
A {\em scaffold} for $L$ is a smooth submanifold $W$ of $M$ with the
following properties:

(1) $\p L \subset W$;

(2) $\bfn\in\Gamma (T_{\p L} W)^\omega$
(here, $S^\omega$ denotes the symplectic orthogonal complement
of a subspace $S$ of a symplectic vector space $V$, defined by
$S^\omega \equiv\{v\in V : \omega(v,s) = 0\ \forall s\in S\})$;

(3) the bundle $(TW)^\omega$ is trivial.
\end{defi}

The deformations of special Lagrangian submanifolds with boundary
constrained to be in a fixed scaffold were studied by Butscher
who proved the following.

\begin{thm}[{\cite[Main Theorem]{butscher}}]\label{bSL}
Let $L$ be a compact special Lagrangian
submanifold of a Calabi-Yau manifold $M$ with non-empty boundary $\p L$ and
let $W$ be a symplectic, codimension two scaffold for $L$. Then the moduli
space of nearby minimal Lagrangian deformations of L with boundary on W is
finite dimensional and is locally parameterised by the vector space of closed
co-closed 1-forms on L satisfying Neumann boundary conditions
$$
{\cal H}^1_{\bfn} (L)=
\{\alpha\in\Omega^1(L):
d\alpha=0,\; d^*\alpha=0,\; (\bfn \lrcorner\alpha)|_{\p L}=0\}.
$$
\end{thm}

Theorem~\ref{bSL} allows special Lagrangian deformations of $L$ with
arbitrary phase $\theta$ and the proof uses an extended version of the
deformation map including $\theta$ as an additional variable.
The deformation map also requires a construction of an auxiliary metric
$\hat{g}$ so that $W$ is totally geodesic for $\hat{g}$ and the appropriate
version of the tubular neighbourhood theorem (Proposition~\ref{tubular}(b)).
(The condition (3) in the definition of the scaffold is used in the
construction of~$\hat{g}$.)
The argument then proceeds by appealing to the implicit function theorem
in a similar manner to the deformation problem for closed submanifolds. 
This requires a version of Hodge theory for compact manifolds with boundary
\cite{Schwarz} to identify appropriate Banach subspaces of forms
so that the linearization of the deformation map be surjective.

The vector space ${\cal H}^1_{\bfn}(L)$ in Theorem~\ref{bSL}
is naturally isomorphic to the real cohomology group $H^1(L,\RE)$ and thus
has dimension $b^1(L)$ by Hodge theory for manifolds with boundary,
see \cite[p.~927]{CDGM}.

Observe that the condition (2) in the definition of a scaffold
means that $J\bfn$ is perpendicular to $W$. This will be
automatically satisfied when $W$ is a complex submanifold of positive
codimension in $M$, so the tangent spaces of $W$ are invariant under~$J$.
In this case, we also have $\Omega|_W=0$. Then, applying
Proposition~\ref{calibr-min}(b), we obtain that for each special
Lagrangian $L$ with boundary in~$W$ the minimal Lagrangian deformations of
$L$ with boundary confined to $W$ will actually be special Lagrangian.
We thus obtain a variant of Butscher's result for the space of special
Lagrangian deformations.

\begin{cor}\label{SLdeform}
Let $L$ be a compact special Lagrangian
submanifold of a Calabi-Yau manifold M with non-empty boundary $\p L$.
Let $W$ be a complex codimension one submanifold of~$M$ with trivial normal
bundle and with $\p L \subset W$ (in particular, $W$ is a scaffold for $L$).
Then the moduli space of nearby special Lagrangian deformations of $L$ is a
smooth manifold of dimension $b^1(L)$.
\end{cor}

\section{Coassociative submanifolds in $G_2$-manifolds}

The two calibrations considered in this and the
next section are defined on 7-dimensional manifolds with a torsion-free
$G_2$-structure. We shall first briefly recall some key definitions.
The readers are referred to \cite{HL}, \cite[Ch.~11,12]{joyce} and the article
by Karigiannis in this volume for a more detailed account of $G_2$-structures
and the related calibrations.

The group $G_2$ can be defined, following
\cite{bryant}, as the stabilizer, in the standard action of $GL(7,\RE)$ on
$\Lambda^3(\RE^7)^*$, of the 3-form
\begin{equation}\label{phi0}
\varphi_0 = dx^{123} + dx^{145} + dx^{167} + dx^{246}
 - dx^{257} - dx^{347} - dx^{356},
\end{equation}
where d$x^{123}=dx^1\we dx^2\we dx^3$ and so on, with $x^1,\ldots,x^7$ the
usual coordinates on the Euclidean~$\RE^7$.
The 3-form in~\eqref{phi0} encodes the cross-product defined by considering
$\RE^7$ as pure imaginary octonions and setting
$\varphi_0(a,b,c) = \langle a\times b, c\rangle$.
The group $G_2$ is a
$14$-dimensional Lie group and a subgroup of $SO(7)$.

The Hodge dual of $\varphi_0$ is a 4-form given by:
\begin{equation*}
*\varphi_0 = dx_{4567}+dx_{2367}+dx_{2345}+dx_{1357}-dx_{1346}-dx_{1256}-dx_{1247}.
\end{equation*}

Let $M$ be a 7-dimensional manifold. We say that a differential 3-form
$\varphi$ on~$M$ is {\em positive}, or is a {\em $G_2$ 3-form}, if for each
$p\in M$ there is a linear isomorphism $\iota_p:\RE^7\to T_pM$ with
$\iota_p^*(\varphi(p))=\varphi_{0}$, where $\varphi_0$ is given in~\eqref{phi0}.
Every $G_2$-structure on~$M$ can be induced by a positive 3-form $\varphi$ and
we shall, slightly informally, say that $\varphi$ {\em is} a $G_2$-structure.
As $G_2\subset SO(7)$, every $G_2$-structure $\varphi$ induces on $M$ a metric
$g(\varphi)$ and orientation and thus also a Hodge star~$*_\varphi$.

The intrinsic torsion of a $G_2$-structure $\varphi$ on $M$ vanishes precisely
when $d\varphi=0$ and $d {*}_\varphi \varphi=0$ (\cite{fernandez-gray}). In
this case, we call $(M,\varphi)$ a {\em $G_2$-manifold}.

The 4-form $*_\varphi \varphi$ defines on each $G_2$-manifold $M$ a calibration
(cf. \cite[\S IV.1.B]{HL}).
In fact, the results in this section only require that the $G_2$ 3-form be
closed $d\varphi=0$. We say that an oriented 4-dimensional submanifold
$X\subset M$ is a {\em coassociative submanifold} if the equality
$*_\varphi\varphi|_X=\vol_X$ is attained. If in addition
$d {*}_\varphi \varphi=0$ holds, then $X$ is calibrated by $*_\varphi \varphi$
and we call $X$ a {\em coassociative calibrated submanifold}.

The following equivalent definition of coassociative submanifolds will be
useful.

\begin{prop}[cf.~\mbox{\cite[Cor.~IV.1.20]{HL}}]\label{ca-defi}
For an orientable 4-dimensional submanifold~$X$ of a \mbox{7-manifold} $M$
with a $G_2$-structure $\varphi\in\Omega^3_+(M)$, the equality
$*_\varphi\varphi|_X=\vol_X$ holds for some
orientation of~$X$ if and only if \mbox{$\varphi|_X=0$.}
\end{prop}

Let $\varphi$ be a closed $G_2$ 3-form on~$M$ and let a submanifold
$X\subset M$ be coassociative. Then the infinitesimal deformations of~$X$ can
equivalently be given by self-dual 2-forms on~$X$ via an isometry of vector
bundles (cf.~\cite[Proposition 4.2]{mclean})
\begin{equation}\label{j-map}
\jmath_X:\mathbf{v}\in N_{X/M}\to
(\mathbf{v}\,\lrcorner\,\varphi)|_X\in\Lambda^2_+T^*X;
\end{equation}
where $\Lambda^2_+T^*X$ denotes the bundle of self-dual 2-forms. (The
corresponding statements in~\cite{mclean} use anti-self-dual 2-forms because
McLean uses a different sign convention for the $G_2$ 3-form.) The map
\begin{equation}\label{mapF}
F:\alpha\in\Omega^2_+(X)\to  \exp_{\bfv}^*(\varphi)\in\Omega^3(X),
\qquad  \bfv=\jmath_X^{-1}(\alpha),
\end{equation}
is defined for `small' $\alpha$, and $F(\alpha)=0$ precisely if
$\exp_{\bf{v}}(X)$ is a coassociative deformation.

The next theorem summarizes the results obtained by McLean about the
deformations of closed coassociative submanifolds. 

\begin{thm}[\mbox{cf.~\cite[Thm. 4.5]{mclean}, \cite[Thm. 2.5]{joyce-salur}}]\label{mclean}
Let $M$ be a 7-manifold with a closed $G_2$-structure $\varphi$
and $X\subset M$ a coassociative submanifold (not necessarily closed). 
\begin{itemize}
\item[\emph{(a)}] Then for each $\alpha\in\Omega^2_+(X)$, one has
$dF|_0(\alpha)=d\alpha$ and the 3-form 
$F(\alpha)$ (if defined) is exact.
\item[\emph{(b)}] If, in addition, $X$ is compact and without boundary then every closed
self-dual 2-form $\alpha$ on~$X$ arises as $\alpha=\jmath_X(\bfv)$, for some
normal vector field~$\bfv$ tangent to a smooth 1-parameter family of
coassociative submanifolds containing~$X$.  Thus, in this case, the space of
nearby coassociative deformations of~$X$ is a smooth manifold parameterized by
the space $\HH^2_+(X)$ of closed self-dual 2-forms on~$X$.
\end{itemize}
\end{thm}
\begin{remark}
Self-dual 2-forms on a compact manifold without boundary are closed precisely
if they are harmonic. By Hodge theory, the dimension of $\HH^2_+(X)$ is
therefore equal to the dimension $b^2_+(X)$ of a maximal positive subspace for
the intersection form on~$X$. It is thus a topological invariant.
\end{remark}

The hypotheses of Theorem~\ref{mclean} do not include the co-closed
condition $d{*}_\varphi\varphi=0$. In fact, the argument in~\cite{mclean}
constructs a smooth moduli space of closed submanifolds $X$ satisfying
$\varphi|_X=0$, for a closed $G_2$-structure $\varphi$. 

There is a certain analogy between McLean's deformation theory of closed
coassociative submanifolds of $G_2$-manifolds and closed special Lagrangian
submanifolds of Calabi--Yau manifolds. In both cases, the respective
submanifolds are calibrated and minimal and have an equivalent definition
in terms of the vanishing of appropriate real differential forms on the
ambient manifold. The deformation theory is `unobstructed' and there is a
smooth finite-dimensional moduli space, locally parameterised by some
finite-dimensional space of harmonic forms on the submanifold with the
dimension a topological invariant obtained by Hodge theory.

On the other hand, when the submanifold has a boundary the deformation
theories become rather different. As we explain below, following~\cite{KL},
the deformation problem for compact coassociative submanifolds with
boundary cannot possibly be set as a boundary value problem of first order
with standard Dirichlet or Neumann boundary conditions. Instead the
deformation problem will be `embedded' in a boundary value elliptic problem
of second order.

A suitable choice of boundary considerations is again facilitated by the
concept of (a coassociative version of) a scaffold which we now define.
By way of preparation, we consider an orientable 6-dimensional submanifold
$S$ in a 7-manifold $M$ with a closed $G_2$ 3-form $\varphi$ on~$M$. 
The normal bundle of $S$ is trivial and there is a `tubular neighbourhood'
$T_S$ of~$S$ diffeomorphic to $S\times\{-\eps<s<\eps\}$, such that $S$
corresponds to $\{s=0\}$ and $\bfn=\frac\p{\p s}$ is a unit vector field
on~$T_S$ with $\bfn|_S$ orthogonal to~$S$ in the metric $g(\varphi)$.
More precisely, we consider a coassociative submanifold with
(compact) boundary contained in $S$ and then the required $T_S$ exists
after shrinking $S$ to some neighbourhood of this boundary.
We can write
$$
\varphi|_{T_S}=\omega_s\we ds + \Upsilon_s,
$$
for some 1-parameter families of 2-forms $\omega_s$ and 3-forms
$\Upsilon_s$ on~$S$.

The forms $\omega_0=(\bfn_S\lrcorner\varphi)|_S$ and $\Upsilon_0=\varphi|_S$
together define an $SU(3)$-structure on~$S$, in general with torsion.
This can be seen point-wise, by a consideration similar to
$G_2$-structures earlier in this section, from the property that the
simultaneous stabilizer of $\omega_0$ and $\Upsilon_0$ in the action of
$GL(T_p S)$ at each $p\in S$ is isomorphic to $SU(3)$. In particular,
$\omega_0^3$ defines an orientation on~$S$.

\begin{defi}
An orientable 6-dimensional submanifold $S$ is a {\em symplectic submanifold}
of a 7-manifold $(M,\varphi)$ with a closed $G_2$-structure if
$d_S\omega_0=0$, where $\omega_0$ is as defined above and $d_S$ is the
exterior derivative on~$S$.

A 3-dimensional submanifold $L\subset S$ of a symplectic submanifold
$S\subset M$ is said to be \emph{special Lagrangian} if 
$\omega_0|_L=0$ and $\varphi|_L=0$.
\end{defi}

Every $S$ in the above definition has a Hermitian symplectic structure
compatible with the $SU(3)$-structure induced from~$M$. In particular, a
non-vanishing $(3,0)$-form on~$S$ is obtained as $\Omega_0=*_S\Upsilon -
i\Upsilon$, where the 6-dimensional Hodge start $*_S$ is taken with respect
to the orientation $\omega_0^3$ and the induced metric from~$M$. As noted
in the previous section, McLean's theory remains 
valid and any closed special Lagrangian $L\subset S$ has a smooth moduli
space of dimension $b^1(L)$ of nearby special Lagrangian deformations.

\begin{defi}\label{ca-scaf}
Let $(M,\varphi)$ be a 7-manifold with a closed $G_2$-structure
and $X\subset M$ a coassociative submanifold with boundary $\p X$.
We say that an orientable 6-dimensional submanifold $S$ of $M$ is a
\emph{scaffold}  
for $X$ if
\begin{itemize}
\item[(a)] $X$ meets $S$ orthogonally, i.e.\
$\p X\subset S$ and the normal vectors to~$S$ at $\p X$
are tangent to~$X$,  $\bfn\in N_{S/M}|_{\p X}$, and
\item[(b)] $S$ is a symplectic submanifold of~$(M,\varphi)$.
\end{itemize}
\end{defi}

One notable property of a scaffold $S$ in Definition~\ref{ca-scaf} is that 
for each coassociative $X$ meeting $S$ orthogonally, the intersection
$L=X\cap S$ is special Lagrangian in~$S$.

The infinitesimal deformations of compact coassociative $X$ with
boundary in a fixed submanifold $S$ correspond via \eqref{j-map} to a
subspace of self-dual 2-forms on~$X$ satisfying boundary conditions.
We can write the restriction of any $2$-form $\alpha$ on $X$ to a collar
neighbourhood $C_{\p X}=T_S\cap X$ of the boundary as
$\tilde\alpha=\alpha_\tau+\alpha_\nu\we ds$.
The Dirichlet and Neumann boundary conditions for $\alpha$ are then given by,
respectively, $\alpha_\tau=0$ and $\alpha_\nu=0$. When $\alpha$ is self-dual,
the two conditions are equivalent and force $\alpha$ and the corresponding
normal vector field $\jmath_X^{-1}(\alpha)$ to vanish at each point of 
$\p X$.\label{eight} However, if $d\alpha=0$ and $\alpha$ vanishes
on the boundary then $\alpha=0$ by \cite[Lemma 2]{CDGM}.
This may be understood as an extension of \cite[Thm. IV.4.3]{HL}, which
states that there is a locally unique coassociative submanifold containing any
real analytic 3-dimensional submanifold upon which $\varphi$ vanishes.

It turns our that a suitable choice of the infinitesimal deformations with
boundary condition is given by the following subspace of self-dual 2-forms
\begin{equation}\label{hsdbc}
\Omega^2_+(X)_{\text{bc}}
=\{\alpha\in 
\Omega^2_+(X):\;
\bfn\lrcorner d\alpha=0\text{ and }d_{\p X}(\bfn\lrcorner\alpha)=0
\text{ on $\p X$}\}.
\end{equation}
The boundary condition $\bfn\lrcorner d\alpha=0$ ensures that every
harmonic form in $\Omega^2_+(X)_{\text{bc}}$ is closed and the 
boundary condition $d_{\p X}(\bfn\lrcorner\alpha)=d_{\p X}\alpha_\nu=0$
means that the boundary $\p X$ will only move in the space of special
Lagrangians in~$S$.

The subspace of $\Omega^2_+(X)_{\text{bc}}$ of the coassociative
infinitesimal deformations is given by the harmonic (or closed) forms
$(\mathcal{H}^2_+)_{\text{bc}}$. It has a finite dimension $\le b^1(\p X)$
and our next theorem asserts that elements of $(\mathcal{H}^2_+)_{\text{bc}}$
`integrate' to actual coassociative deformations of~$X$.

\begin{thm}[{\cite[Thm. 1.1]{KL}}]\label{thm1}
Suppose that $M$ is a 7-manifold with a $\text{\emph{G}}_2$-structure given
by a closed 3-form. The moduli space of compact coassociative local deformations
of~$X$ in $M$ with boundary $\p X$ in a scaffold $S$ is a
finite-dimensional smooth manifold parameterized by
$(\mathcal{H}^2_+)_{\text{bc}}$. The dimension of this moduli space is not
greater than~$b^1(\p X)$.
\end{thm}

Here is an example when strict inequality
$\dim (\mathcal{H}^2_+)_{\text{bc}} < b^1(\p X)$ occurs.

\begin{example}[{\cite[p. 72]{KL}}]
A K\"ahler complex 3-fold $(Z,\omega)$ is called \emph{almost Calabi--Yau} if
it admits a nowhere vanishing holomorphic $(3,0)$-form~$\Omega$. 
Then the 7-manifold $M=Z\times S^1$ has a closed
$G_2$-structure $\omega\we d\theta+\mathrm{Re}\,\Omega$, where $\theta$
is a coordinate on~$S^1$. Let $X=L\times S^1\subset M$ be a
compact coassociative 4-fold. Then $L$ is special Lagrangian in $Z$.  
We can think of $X$ as an embedding of a manifold $L\times[0,1]$ whose two
boundary components, $L\times\{0\}$ and $L\times\{1\}$, are mapped to
$L$ in $Z$.
(If $Z$ is Calabi--Yau, i.e.\ if \eqref{omegas} also holds, then $X$ is
a calibrated coassociative.)
It is not difficult to see that $Z\times\mathrm{pt}$ is a
scaffold for $X$. Theorem \ref{thm1} gives us that $X$ has a smooth
moduli space of coassociative deformations with
dimension $\le 2b^1(L)$.

Let $\alpha\in\Omega^2_+(X)$.  Then $\alpha=\xi_\theta\we d\theta+*_L\xi_\theta$,
for some path of 1-forms $\xi_\theta$ on~$L$.  
It follows from \cite[Thm. 3.4.10]{Schwarz} that a harmonic self-dual 2-form
on $X$ is uniquely determined by its values $\xi_0$, $\xi_1$ on the boundary.
The subspace of harmonic $\alpha\in\Omega^2_+(X)$ such that $\xi_0$ and $\xi_1$
are harmonic on~$L$ has dimension $2b^1(L)$ and corresponds
precisely to the paths $\xi_\theta=(1-\theta)\xi_0+\theta\xi_1$.
On the other hand, $\alpha\in (\HH^2_+)_{\text{bc}}$ if and only if $\alpha$
is harmonic and ${\partial\xi_{\theta}}/{\partial\theta}=0$, 
so $\xi_0=\xi_1$. Thus
$\dim(\HH^2_+)_{\text{bc}}=b^1(L)<b^1\big((L\times\{0\})\sqcup 
(L\times\{1\})\big)$ in this example.

This can also be seen geometrically.  If the deformations of the aforementioned 
two boundary components coincide in~$Z\times\mathrm{pt}$ then, by taking a
product with $S^1$, we obtain a coassociative deformation of
$X=L\times S^1$ defining a point in the moduli space in Theorem~\ref{thm1}.
On the other hand, if a coassociative deformation $\tilde{X}$ of $X$ is such
that the deformations $\tilde{L}_0$ and $\tilde{L}_1$ of $L\times\{0\}$ and
$L\times\{1\}$ are special Lagrangian but {\em distinct} then  $\tilde{X}$ and
$\tilde{L}_0\times S^1$ are two distinct coassociative 4-folds intersecting in
a real analytic 3-fold on which $\varphi$ vanishes, which 
contradicts \cite[Thm. IV.4.3]{HL}. Therefore, the moduli space in this example
is identified with special Lagrangian deformations of $L$ in the almost
Calabi--Yau manifold $Z$. As we noted earlier, these deformations
have a smooth moduli space of dimension~$b^1(L)$.
\end{example}

\section{Associative submanifolds in $G_2$-manifolds}

Let $M$ be again a 7-dimensional manifold with a $G_2$-structure given by a
positive 3-form $\varphi$, as defined in the previous section. If $d\varphi=0$
then $\varphi$ is a calibration on~$M$ as $\varphi_p|_V\le \vol_V$ for each
oriented 3-plane in $T_p M$ (cf. \cite[\S IV.1.A]{HL}). The equality is attained
precisely when $V$ is an {\em associative subspace} of $T_p M$, i.e.\ is
closed under the cross-product on $T_p M$ induced by the $G_2$-structure
$\varphi$. (The corresponding subalgebra $(V,\times)$ is isomorphic to $\RE^3$
with the standard vector product.)

The deformation theory discussed in this section does not always require the
$G_2$ form $\varphi$ to be closed. Similarly to the discussion of
coassociative submanifolds in the previous section, we shall define the term
{\em associative submanifold} $Y$ in $(M,\varphi)$ for an arbitrary $G_2$
structure $\varphi$, meaning a 3-dimensional submanifold $Y$ satisfying
$\varphi|_Y=\vol Y$. If also the $G_2$-structure is closed $d\varphi=0$, then
we shall call $Y$ an {\em associative calibrated submanifold}; indeed, in this
case $Y$ is calibrated by~$\varphi$.

McLean~\cite{mclean} studied the deformation theory of closed associative
calibrated submanifolds in $G_2$-manifolds. His results were later
extended by Akbulut and Salur \cite{AS} to arbitrary $G_2$-structures on
7-manifolds. We assume a torsion-free $G_2$-structure to simplify some
details.

As in the previous sections, it is useful to first note that associative
submanifolds of $M$ can be equivalently defined by the vanishing of an
appropriate differential form on~$M$. In the present case, it is a 3-form
$\chi$ with values in $TM$,
$$
\chi = \sum_{k=1}^7 (\eta_k\lrcorner *_\varphi\varphi)\otimes\eta_k,
$$
for any local orthonormal positively oriented frame field
$(\eta_k)_{k=1}^7$ on~$M$  \cite[p. 1217]{gayet}.

\begin{prop}[{cf. \cite[\S 5]{mclean}}]
A 3-dimensional submanifold $Y$, with some choice of
orientation, is associative if and only if $\chi|_Y=0$.
\end{prop}

Suppose that $Y$ is a closed associative submanifold in a $G_2$-manifold
$(M,\varphi)$ and let $\bfv\in \Gamma(N_{Y/M})$ be a normal vector field
along~$Y$. The deformation map for $Y$ is defined as
$$
F: \bfv \in\Gamma(N_{Y/M}) \to \exp_{\bfv}^*\tau\in\Omega^3(Y,TM|_Y)
$$
The linearization of $F$ at $\bfv=0$ is given by
\begin{equation}\label{dirac}
D\bfv = \sum_{i=1}^3 e_i \times \nabla_{e_i}^\bot \bfv,
\end{equation}
where $e_1,e_2,e_3$ is any positively oriented local orthonormal frame
field of~$TY$ (thus $e_3=e_1\times e_2$), the cross-product is induced by
$\varphi$ and a connection $\nabla^\bot$ on
$N_{Y/M}$ is induced by the Levi--Civita connection of $(M,g(\varphi))$.
Since $Y$ is 3-dimensional and associative, both $TY$ and $N_{Y/M}$ are
trivial vector bundles \cite[Remark 2.14]{joyce-karigiannis} and the 
expression~\eqref{dirac} is valid globally over~$Y$. There is an
invariant interpretation of $N_{Y/M}$ as a vector bundle associated with a
principal $\Spin(4)$-bundle over $Y$ via the tensor product of a spin
representation and some other representation. Then $D$ becomes the respective
Dirac type operator (meaning that the principal symbol of $D^2$ is
$\sigma(D^2)(p,\xi)=\|\xi\|^2$, for all $p\in Y$).

The map $F$ makes sense for an arbitrary $G_2$-structure $\varphi$ on~$M$,
but when the $G_2$-structure is not torsion-free the expression 
\eqref{dirac} for~$D$ then has extra terms of order zero. So in this
more general case $D$ is still a Dirac type operator with the same
principal symbol.

We thus obtain.
\begin{thm}[{\cite{mclean,AS}}]
For a closed associative submanifold $Y$ in a $G_2$-manifold $(M,\varphi)$,
the Zariski tangent space to associative deformations
of~$Y$ is finite-dimensional, given by the kernel of the Dirac type
operator $D$ in~\eqref{dirac}, an elliptic operator of index~$0$. In
particular, $Y$ is either rigid or the associative deformations of $Y$ are
obstructed, i.e.\ a section $\bfv$ with $D\bfv=0$ need not arise as
$\dot{s}_0$ from any 1-parameter family $Y_t=\exp s_t$ of associative
submanifolds with $s_0=0$.
\end{thm}

The deformations of compact associative submanifolds with boundary
contained in a fixed submanifold (scaffold) were investigated by Gayet and
Witt \cite{GW}, see also Gayet~\cite{gayet}. An appropriate choice of scaffold
in this situation is given by a 4-dimensional submanifold $X$ such that no
tangent space $T_p X$ contains associative 3-planes. In particular, $X$ may
be any coassociative submanifold and we shall assume this below for technical
convenience.
In this case, any associative calibrated submanifold with boundary in~$X$ is
volume minimizing in its relative homology class as $\varphi|_X=0$ by
Propositions \ref{calibr-min}(b) and \ref{ca-defi}.

Let $Y$ be a compact associative submanifold with boundary $\p Y$ contained
in a coassociative submanifold~$X$. Denote by $\bfn\in\Gamma(TY|_{\p Y})$
the inward-pointing unit normal along $\p Y$.
For each point $p\in \p Y$ the cross-product of the $G_2$-structure $\varphi$
$$
J(v)=\bfn_p\times v
$$
defines an (orthogonal) complex structure on the orthogonal complement
$(\bfn_p)^\bot\subset T_pM$. Further, $J$ acts on the fibres of the normal
bundle $N_{\p Y/X}$, making it into a complex line bundle.
Note that $N_{\p Y/X}$ is a subbundle of $N_{Y/M}|_{\p Y}$; the respective
orthogonal complement $\mu_{\p Y}$ is also invariant under $J$ and can be
considered as a complex line bundle.
The tangent spaces of $\p Y$ are also preserved by $J$ and in this way
$\p Y$ is made into a compact Riemann surface. 
The latter complex line bundles satisfy an adjunction-type relation
$\bar{\mu}_{\p Y}\cong N_{\p Y/X}\otimes_{\CX} T\p Y$  \cite[Lemma 3.2]{GW}.

The infinitesimal associative deformation problem for $Y$ with boundary
confined to $X$ can be expressed as
\begin{equation}\label{dirac3}
D\bfv = 0, \quad B(\bfv|_{\p Y})=0,\qquad \bfv\in\Gamma(N_{Y/M}),
\end{equation}
where $D$ is the Dirac type operator in~\eqref{dirac} and the zero order
operator $B$ is induced by the orthogonal projection $N_{Y/M}|_{\p Y}\to
\mu_{\p Y}$ with kernel $N_{\p Y/X}$. 

Gayet and Witt prove.
\begin{thm}[{\cite[Thm. 4.4, Cor. 4.5]{GW}}]\label{GW.ind}
Let $(M,\varphi)$ be a 7-manifold endowed with a $G_2$-structure and
let $Y\subset M$ be a compact associative submanifold with boundary contained
in a coassociative submanifold $Y\subset M$.

Then the linear operator $D\oplus B$ in \eqref{dirac3}
defines an elliptic boundary value problem with finite Fredholm index
\begin{equation}\label{index3}
\Ind (D\oplus B) = \sum_j ( \int_{\Sigma_j} c_1(N_{\Sigma_j/X}) + 1 - g_j ).
\end{equation}
Here $\Sigma_j$ denote the boundary components of $\p Y$ with $g_j$ the
genus of $\Sigma_j$, and $c_1(N_{\Sigma_j/X})$ is the first Chern class of the
complex line bundle $N_{\Sigma_j/X}$.
\end{thm}
The key point in the proof of Theorem~\ref{GW.ind} is that the index in
question can be computed as the index of the Cauchy--Riemann operator
$\bar\p_{\p Y/X}$ associated with the complex structure $J$ on $N_{\p Y/X}$.

\begin{example}[{\cite[p.2364]{GW}}]
When the associative submanifold has a boundary, the index in
Theorem~\ref{GW.ind} can be positive. One simple example uses a construction
by Bryant and Salamon~\cite{BS} of torsion-free $G_2$-structure $\varphi$
inducing a complete metric with holonomy $G_2$ on $\mathcal{S}=
S^3\times\RE^4$, the total
space of the spinor bundle over the standard round 3-sphere. The zero section
$S^3\times \{0\}$ is an associative submanifold, being the fixed locus of the
$G_2$-involution acting as $-1$ on the fibres
(note~\cite[Prop. 12.3.7]{joyce}).  
Take $Y\subset S^3\times \{0\}$ to be a 3-dimensional ball, so $\p Y=S^2$.
Let $a$ be a nowhere vanishing section of $\mathcal{S}|_{\p Y} =
\p Y\times\RE^4\to \p Y$. Then $a,\; Ja$ (with $J$ as defined above) generate
a trivial complex line bundle, denote its total space by $\tilde{X}$. It can
be checked that there is a local coassociative submanifold
$X\subset\mathcal{S}$ containing $\p Y$ and with $T_p X=T_p\tilde{X}$ at each
$p\in\p Y$. Theorem~\ref{GW.ind} then applies and the deformation problem
has index~1. This example generalizes, under additional assumptions, to a
complex line bundle $\tilde{X}$ having positive degree~$n$, then the
respective index is $n+1$ \cite[p.2364]{GW}.
\end{example}

Gayet and Witt also gave a generalization of Theorem \ref{GW.ind} where
$X$ is only required to contain no associative 3-planes in its tangent spaces.
The latter property is preserved under small perturbations of the
$G_2$-structure. The deformation problem remains elliptic and the index
formula \eqref{index3} still holds (with appropriate modification of the
definition of the boundary operator~$B$).
This does not guarantee a smooth moduli space of associative
deformations even when the index of $D\oplus B$ is non-negative.
However, Gayet \cite[Thm. 1.4]{gayet} proved that a smooth moduli space of
dimension $\Ind{D\oplus B}$ can be obtained by arbitrary small generic
perturbation of the scaffold~$X$.

\section{Cayley submanifolds in $\Spin(7)$-manifolds}

The calibration considered in this section is defined on 8-dimensional
manifolds with a torsion-free $\Spin(7)$-structure. We begin with a short
summary of the $\Spin(7)$-structure and the Cayley calibration
and refer to \cite{HL} and \cite[Ch.~11,12]{joyce} for further details.
There is a certain, though only partial, analogy with the `geometries'
considered in the previous sections.

The group $\Spin(7)$ can be defined, following \cite{bryant}, as the
stabilizer, in the standard action of $GL(8,\RE)$ on $\Lambda^4(\RE^8)^*$,
of the 4-form $\Phi_0$ written in standard coordinates as
\begin{multline}\label{std.spin7}
\Phi_0=dx_{1234}+dx_{1256}+dx_{1278}+dx_{1357}-dx_{1368}
-dx_{1458}-dx_{1467}\\
-dx_{2358}-dx_{2367}-dx_{2457}+dx_{2468}+dx_{3456}+dx_{3478}+dx_{5678},
\end{multline}
where d$x^{1234}=dx^1\we dx^2\we dx^3\we dx^4$ and so on.
The form $\Phi_0$ arises by considering $\RE^8$ as the (normed) algebra
of octonions, or Cayley numbers, and setting
$\Phi_0(x,y,z,w)=\frac12\langle x(\bar{y}z)-z(\bar{y}x),w\rangle$.
In this way, $\Spin(7)$ is also identified as a subgroup of $SO(8)$.
This form is also self-dual $*\Phi_0=\Phi_0$ with respect to the standard
Euclidean metric and orientation.

Given an oriented 8-manifold~$M$, define a subbundle of
4-forms $\AM M\subset \Lambda^4 T^* M$ 
with the fibre $\AM_p M$ at each $p\in M$ being the set of all 4-forms that
can be identified with $\Phi_0$ via an orientation-preserving isomorphism
$T_pM\to\RE^8$. The fibres of $\AM M$ are diffeomorphic to the orbit
$GL_+(8,\RE)/\Spin(7)$ of $\Phi_0$, a 43-dimensional submanifold of
the 70-dimensional vector space $\Lambda^4 (\RE^8)^*$.

A choice of 4-form $\Phi\in\Gamma(\AM M)$ is equivalent to a choice
of a $\Spin(7)$-structure on $M$. By the inclusion $\Spin(7)\subset SO(8)$,
every such $\Phi$ induces on $M$ a metric $g=g(\Phi)$, an orientation and a
Hodge star $*_\Phi$, with ${*_\Phi}\Phi=\Phi$. We shall sometimes
refer to $\Phi$ as a $\Spin(7)$-structure.

When a form $\Phi\in \AM M$ is closed, $d\Phi=0$, we say that the
$\Spin(7)$-structure $\Phi$ is torsion-free and that $(M,\Phi)$ is a
{\em $\Spin(7)$-manifold}. The condition $d\Phi=0$ is equivalent to the
metric $g(\Phi)$ being Ricci-flat with reduced holonomy contained in $\Spin(7)$
\cite{fernandez}. 
In this case, $\Phi$ defines a calibration on~$M$.

Given a $\Spin(7)$-structure on an 8-manifold $M$, we say that an oriented
4-dimensional submanifold $P\subset M$ is a {\em Cayley submanifold} if
$\Phi|_P=\vol_P$ We say that $P$ is a {\em Cayley calibrated submanifold} if
in addition $d\Phi=0$, i.e.\ precisely if $P$ is calibrated by $\Phi$.

Let $P$ is a compact Cayley calibrated submanifold of a $\Spin(7)$-manifold
$(M,\Phi)$. The infinitesimal Cayley deformations of $P$ are given by the
kernel of a first order elliptic operator $D:\Gamma(N_{P/M})\to \Gamma(E)$,
for a rank 4 vector bundle $E$ over~$P$,
$$
E = \{\alpha\in\Lambda^2_7 M|_P : \alpha|_{TP}=0\}.
$$
Here $\Lambda^2_7 M\subset \Lambda^2 T^* M$ is a subbundle corresponding to
an irreducible representation of $\Spin(7)$ on the space of 2-forms
$\Lambda^2(\RE^8)^*$.
The following result was proved by McLean.
\begin{thm}[{\cite[Thm. 6.3]{mclean}}]\label{zariski8}
Let $P$ be a Cayley submanifold of a $\Spin(7)$-manifold $(M,\Phi)$.
Then the Zariski tangent space to Cayley deformations
of~$P$ is finite-dimensional, given by the kernel of the elliptic operator
\begin{equation}\label{dirac8}
D:\bfv\in\Gamma(N_{P/M}) \to
\sum_{i=1}^4 e_i \times \nabla_{e_i}^\bot \bfv \in \Gamma(E),
\end{equation}
where $e_1,e_2,e_3,e_4$ is any positively oriented local orthonormal frame
field of~$TP$, the cross-product is induced by the $\Spin(7)$-structure $\Phi$
and a connection $\nabla^\bot$ on $N_{P/M}$ is induced by the Levi--Civita
connection of $(M,g(\Phi))$.
\end{thm}
\begin{remarks}
When $P$ is a spin manifold, there is an invariant interpretation of $D$ using
a spin structure on~$P$ \cite[\S 6]{mclean}. Denote by $\mathcal{S}_+$ and
$\mathcal{S}_-$ the 
positive and negative spinor bundles over~$P$. Then
\begin{equation}\label{Cayley-spin}
N_{P/M}\otimes_\RE \CX \cong \mathcal{S}_+\otimes_\CX F
\text{ and }
E \otimes_\RE \CX \cong \mathcal{S}_-\otimes_\CX F,
\end{equation}
for some quaternionic line bundle~$F$ over $P$. The operator $D$ is
identified, via~\eqref{Cayley-spin}, with a positive Dirac type operator
associated with a connection on~$F$.

Theorem~\ref{zariski8} was extended to arbitrary $\Spin(7)$-structures in
\cite[\S 13]{GIP}. When the $\Spin(7)$-structure is not torsion-free, the
expression \eqref{dirac8} for~$D$ has extra terms of order zero. This does not
affect the principal symbol or the index of~$D$.
\end{remarks}

Deformations of Cayley submanifolds were further investigated by
Ohst and included the following.
\begin{thm}[{\cite[Prop. 3.4 and Thm. 3.10]{ohst}}]\label{gen-closed}
Let $(M,\Phi)$ be an 8-manifold with a $\Spin(7)$-structure and $P\subset M$ a
closed Cayley submanifold. Then
\begin{itemize}
\item[\emph{(a)}]
 the index of the operator~\eqref{dirac8} associated with $P$ is
\begin{equation}\label{index8}
\Ind D = {\textstyle\frac12}\chi(P) + {\textstyle\frac12}\sigma(P)
- [P]\cdot [P],
\end{equation}
where $\chi(P)$ is the Euler characteristic, $\sigma(P)$ is the signature and 
\mbox{$[P]\cdot [P]$} is the self-intersection number of~$P$.
\item[\emph{(b)}]
For every generic $\Spin(7)$-structure $\tilde{\Phi}$ on $M$ such that
$\|\tilde\Phi - \Phi\|$ is sufficiently small and $\tilde{\Phi}$ induces the
same metric $g(\tilde{\Phi})=g(\Phi)$ the following holds.
The moduli space of Cayley submanifolds with respect to $\tilde{\Phi}$ which
are $C^{1,\alpha}$-close to $P$ ($0<\alpha<1$) is either empty or a smooth
manifold of dimension $\Ind D$ (if $\Ind D \ge 0$).
\end{itemize}
\end{thm}
\begin{remarks}
The submanifold $P$ need not be Cayley with respect to $\tilde\Phi$, thus $P$
need not be in the respective moduli space. 
The norm $\|\tilde\Phi - \Phi\|$
can be taken to be the $C^{1,\alpha}$-norm on a compact neighbourhood of~$P$.

The variant of Theorem~\ref{gen-closed} also holds with $\tilde{\Phi}$ 
generic in the set of all $\Spin(7)$-structures close to $\Phi$, i.e.\ without
the restriction on the metric $g(\tilde{\Phi})$.
\end{remarks}

We next turn to deformations of compact Cayley submanifolds with boundary in a
fixed submanifold (scaffold). The result given below is again due to Ohst
and allows a range of dimensions of the scaffold.

\begin{thm}[{\cite[Thm. 4.18]{ohst}}]\label{b-Cayley}
Let $(M,\Phi)$ be a Spin(7)-manifold and $W$ a submanifold of $M$ with
$3\le\dim W\le 7$. Let $P$ be a compact, connected Cayley submanifold of
$M$ with non-empty boundary $\p P\subseteq W$ such that $P$ and $W$ meet
orthogonally.

Then for every generic torsion-free Spin(7)-structure $\tilde\Phi$ which is
$C^{2,\alpha}$-close to $\Phi$, the moduli space of all Cayley (calibrated)
submanifolds in $(M,\tilde\Phi)$ which are $C^{2,\alpha}$-close to $P$ and
have boundary contained in~$W$ and meet $W$ orthogonally in the metric
$g(\tilde{\Phi})$ is a finite set (possibly empty). Here $0<\alpha<1$.
\end{thm}

The proof of Theorem~\ref{b-Cayley} uses a second order elliptic boundary
problem implied by the linearization $D$ in~\eqref{dirac8} of the Cayley
deformation map. This is because $D$ admits no suitable elliptic boundary
conditions \cite[Prop. 4.21]{ohst}.

The condition $\Phi|_W=0$ required in Proposition~\ref{calibr-min}(b) for the
volume-minimizing property in the relative homology class of $P$ can only hold
if $\dim W\le 4$. Suppose that $\p P$ is a deformation retract of~$W$.
If $\dim W=5$, then $v=d(*_W(\Phi|_W))^\flat$ restricts to a vector field on
$\p P$. If, further, $v$ is parallel with respect to the induced metric
on~$W$, then $P$ is volume-minimizing among the nearby deformations in its
relative homology class.

When $\dim W=6$ and $d(*_W(\Phi|_W))=0$, the scaffold $W$ is a symplectic
submanifold. Then $P$ minimizes the volume among all the submanifolds $P'$
in its relative homology class, with boundary $\p P'\subset W$
a Lagrangian nearby deformation of $\p P$. In all of the above situations,
the minimal volume in the relative homology class of~$P$ is attained
precisely by Cayley calibrated submanifolds \cite[\S 5.1]{ohst}.

\end{document}